\documentclass[12pt, fleqn]{article}
\usepackage[cp1251]{inputenc}
\usepackage{latexsym,amsfonts,amssymb}
\usepackage{graphicx}

\usepackage{amsbsy}
\usepackage{amsmath}
\usepackage{epsf}

\newtheorem{theo}{Theorem}

\newcommand{\bt}{\begin{theo}}
\newcommand{\et}{\end{theo}}
\newcommand{\bd}{\begin{displaymath}}
\newcommand{\ed}{\end{displaymath}}

\newcommand{\be} {\begin{equation}}
\newcommand{\ee} {\end{equation}}
\newcommand{\ba} {\begin{array}}
\newcommand{\ea} {\end{array}}
\newcommand{\bea}{\begin{eqnarray}}
\newcommand{\eea} {\end{eqnarray}}

\begin{document}

\begin{center}
 {\Large \bf Comments on the paper ``Exact solutions of nonlinear diffusion-convection-reaction equation: A Lie symmetry approach'' }
% \vspace{0.5cm}

\medskip

{\bf Roman Cherniha \footnote{\small  Corresponding author. E-mail: r.m.cherniha@gmail.com}}
 \\
{\it ~Institute of Mathematics,  National Academy
of Sciences  of Ukraine,\\
 3, Tereshchenkivs'ka Street, Kyiv 01004, Ukraine
}
 \end{center}

\begin{abstract}
This comment is devoted to the    paper ``Exact solutions of nonlinear diffusion-convection-reaction equation: A Lie symmetry approach'' (CNSNS, {\bf 67} (2019), 253-263) in which several results are not new because  were  derived much earlier. Moreover, some results in the paper   are incorrect or  incomplete.
%%or misleading.}
\end{abstract}


\section{Lie symmetry analysis} \label{sec-1}

The recent paper \cite{molati-mura-2019} is devoted to search for Lie symmetries and  exact solutions of the
diffusion-convection-reaction equation
\be \label {1}
 {u_t}= \ (u^{m})_{xx}+  (b_0u+ b_1u^{p+1})_x
 +(1-u^p)(c_0+c_1u^p)u^{2-m},\ee
 where $u(t,x)$ is unknown smooth function and all the parameters $m,b_0\dots,p$ are real constants (Authors assume that $m>0$ but the sign of $m$  does not play any role  for search of Lie symmetries).

 It can be noted from the very beginning  that the linear term $b_0u$  can be removed from Eq.(\ref{1}) by the well-known transformation, the  Galilei boost
  \be \label{2}
  x^*=x+b_0t.
  \ee
  Hence Eq.(\ref{1})  simplifies to the form
  \be \label {2**}
 {u_t}= \ (u^{m})_{x^*x^*}+  b_1(u^{p+1})_{x^*}
 +(1-u^p)(c_0+c_1u^p)u^{2-m}.\ee
 In Lie symmetry analysis, especially for the Lie symmetry classification, such type  of transformations  plays important role and are called equivalence transformations (ETs)(see, e.g., monographs \cite{ovs80,ch-se-pl-book}). If one ignores ETs then Lie symmetry analysis degenerates in a chaotic process with a lot of  equivalent  cases. Paper \cite{molati-mura-2019} is a typical example.

 To the best of my knowledge, transformation (\ref{2}) was firstly  identified as an ET of the general diffusion-convection-reaction equation (not only Eq.(\ref{1})!)
 \be \label{1*}
u_t=[A(u)u_x]_x+B(u)u_x+C(u), \ee
(here  $A(u)\not=0, B(u)$ and $ C(u)$ are arbitrary smooth functions)
 in paper \cite{ch-se-98}(see formulae (2.12)-(2.13) therein). In \cite{ch-se-ra-08}, the full group of ETs for Eq.(\ref{1*})  was   derived using a rigorous algorithm. The group has the form
 \be\label{2*}
  t^*=\kappa_0t+d_0,\quad  x^*=\kappa_1x+gt+d_1,\quad  u^*= \kappa_2u+d_2,\ee
  where all the parameters $\kappa_0, \dots, d_2$ are arbitrary constants and $\kappa_0\kappa_1\kappa_2 \not=0$.
  Of course,  the group of ETs  contains transformation (\ref{2}) as a particular case (see the parameter $g$ in (\ref{2*})).

 Now we demonstrate that all the results derived in Section 2.2 of \cite{molati-mura-2019} are particular cases of those from \cite{ch-se-ra-08} (actually they can be also deducted  from earlier paper \cite{ch-se-98}).

 Let us take Eq.(24) and Lie symmetries (25) from paper \cite{molati-mura-2019}. Setting $b_0=0$ (because ET(\ref{2}) removes $b_0$) we obtain (star * is omitted for simplicity in what follows)
 \be \label {4}
 {u_t}= \ (u^{m})_{xx}+  b_1(u^{p+1})_x
\ee
 and
 \be \label {5}
  X_1=\partial_{t},  \, X_2=\partial_{x},  \,     X_3=(m-2p-1)t\partial_{t}+(m-p-1)x\partial_{x})-t\partial_{x}+u\partial_{u}
\ee
 respectively. It can be easily checked  that Eq.(\ref{4}) and Lie algebra (\ref{5})  follow as a particular case  from Case 23 of Table 2 \cite{ch-se-ra-08}, when one sets $\lambda_8=0$ and makes renaming $k\to m-1$ and $m \to p$   \cite{ch-se-ra-08}.

 In a quite similar way, it can be checked that equation  and Lie symmetries from (28) \cite{molati-mura-2019}  can be essentially simplified by removing $b_0$. Having this done,  one again sees that Case 24 of Table 2 \cite{ch-se-ra-08} contains the above result. It should be noted that the parameter $c_1$ in (28) \cite{molati-mura-2019} can be reduced to $\pm 1$ by an appropriate ET.

 Equations  and Lie symmetries from (29)  and (31) \cite{molati-mura-2019} after removing $b_0$ also  follow as particular cases from  Case 23 of Table 2 \cite{ch-se-ra-08}, if one sets $m=k$  and $m=k/2$ \cite{ch-se-ra-08}, respectively.

 Finally, equation  and Lie symmetries from  (30) \cite{molati-mura-2019} are known at least  70 years because after removing $b_0$, one arrives at the classical result published  by Ovsiannikov in 1959(see his book \cite{ovs80}). In \cite{ch-se-ra-08}, this result is presented in Case 10 of  Table 2.

 Thus, all the equations and Lie symmetries presented in  \cite{molati-mura-2019}  were discovered many years ago  in the papers \cite{ch-se-98} and \cite{ch-se-ra-08}. Moreover, the Lie symmetry classification of Eq.(\ref{1}) is incomplete. For example, there is a special case, Eq.(\ref{1}) with  $p+1=m$  and $c_0=0$. It can be noted that the equation with the same structure arises in Case 13 of Table 2 \cite{ch-se-ra-08} and admits four-dimensional Lie algebra. This special case is missed in \cite{molati-mura-2019}.

 \section{Reduction and exact solutions} \label{sec-2}

 There are three main techniques for constructing exact solutions by applying Lie symmetries (see the relevant  discussion about their applicability  in Chapter 1 of \cite{ch-se-pl-book}). Many examples of their applications to a wide range of PDEs are presented in monograph \cite{f-s-s}.
 The first technique is rather trivial. One takes the Lie group corresponding to the known Lie symmetry operator and a known solution of the equation in question. Applying the Lie group to the  known solution, one obtains a set of solutions involving at least one free parameter. In the case of Lie algebras of high dimensionality, very nontrivial formulae of multiplication of exact solutions can be derived, see, e.g. \cite{f-s-s}, \cite{fu-ch-89}.

  The second technique is a direct application of a known Lie symmetry for deriving a special substitution, usually called ansatz, which reduces the given two-dimensional PDE to an ODE. If we start from the the general linear combination of all the symmetries of a given PDE   then all possible inequivalent reductions of this PDE to ODEs can be derived. In the case  of Lie algebras of low dimensionality, it is a  simple task (several examples can be found in \cite{ch-se-pl-book}). However, the task transforms  into a difficult  problem if the equation in question admits Lie algebra of high dimensionality (say, 6 or higher) with a nontrivial structure.

  The third  technique is the most sophisticated one and was used in \cite{molati-mura-2019}. This technique is based on the so-called optimal systems of one-dimensional subalgebras and was suggested by Ovsiannikov \cite{ovs80}. Note that one needs to construct also optimal systems of two-  and higher-dimensional  subalgebras in the case of multidimensional PDEs. In the case of Lie algebras of low dimensionality, these optimal systems are known and were summarized in the classical paper \cite{pat-win-77}(incidentally this paper is not cited in   \cite{molati-mura-2019}), while problems occurring  for Lie algebras of high dimensionality is extensively discussed, e.g., in the recent book \cite{fedorchuky-2018}.
  So, if the Lie algebra of invariance is two-, three-, or  four-dimensional then it is enough to identify this algebra in  \cite{pat-win-77}  and the relevant optimal system of one-dimensional subalgebras can be readily written down. Notably, the authors of \cite{pat-win-77} use the notion 'system of non-conjugated subalgebras' instead of Ovsiannikov's terminology 'optimal system of subalgebras'.

  Although all the Lie algebras obtained in \cite{molati-mura-2019} are 3- or 4-dimensional, Authors decided to re-discover the optimal systems of one-dimensional subalgebras. Unfortunately, the results presented in Table 3 are incorrect, excepting Lie algebra (25). It may be noted that this algebra up to a constant multiplier coincides with one $A^a_{35}$ in Table I \cite{pat-win-77}  and the relevant optimal systems of one-dimensional subalgebras do coincide.
  However, the optimal systems  for Lie algebras (28), (29) and (31) are obviously  wrong because they consist of 5  subalgebras. It can be easily seen from the last column of  Table I \cite{pat-win-77} that the optimal system of any  3-dimensional Lie algebra consists of 4 (not 5!) subalgebras at maximum. Similarly,  the optimal system of one-dimensional subalgebras for Lie algebra (30) is also incorrect because one consists of 8 subalgebras while one consists  of 7  (see the algebra $2A_{2}$ in Table II \cite{pat-win-77}). Moreover, the correct optimal system should contain a subalgebra with an arbitrary parameter (not only $\delta=\pm1$ !).

  Because the optimal systems derived in Table 3 \cite{molati-mura-2019} are incorrect (excepting the first line), the number of reduced equations listed on Page 12 \cite{molati-mura-2019} should be essentially smaller. Indeed, according to the general theory (see for details \cite{ovs80}), many solutions of those ODEs will lead to the same solutions of the corresponding nonlinear PDEs.

  Finally, it should be noted that all the exact solutions obtained in Section 3.2  were known earlier. If one sets  $b_0=0$ in the solutions derived then it can be easily identified. For example, the most complicated  exact solution (47) with $b_0=0$ is nothing else but a particular case of that obtained in  \cite{ch-pl-07}(see formula (71) therein).

\section{Conclusions} \label{sec-3}

In this comment, it is shown that all the nonlinear PDEs and their Lie symmetries derived in the recent paper \cite{molati-mura-2019} are  equivalent to those derived earlier in \cite{ch-se-ra-08, ch-se-98} (the detailed proofs and several applications are presented in the recent monograph \cite{ch-se-pl-book}). Moreover,  the optimal systems of one-dimensional subalgebras of Lie algebras are incorrect and do not coincide (excepting a single case) with the results of the seminal work \cite{pat-win-77}. The exact solutions obtained in \cite{molati-mura-2019} also follows from earlier works.

In conclusion, I would like to stress that nowadays  there are many papers  devoted to the Lie symmetry analysis  of nonlinear PDEs, in which the authors simplify the analysis to primitive calculations without knowing state-of-art (see, e.g.,  the recent books \cite{ch-se-pl-book, fedorchuky-2018,ch-dav-book,bl-anco-10} devoted to the symmetry-based methods  and its direct  applications). As a result, 'new' Lie symmetries and 'new' exact solutions are either equivalent to derived in earlier papers, or simply wrong,  optimal systems of subalgebras   are incorrect and real-world  applications are absent (see, e.g., another  typical example commented  in \cite{ch-2020}).

\end{document}